# FLAT VECTOR BUNDLES OVER PARALLELIZABLE MANIFOLDS

JÖRG WINKELMANN

ABSTRACT. We study flat vector bundles over complex parallelizable manifolds.

## 1. INTRODUCTION

Let $G$ be a complex Lie group and $\Gamma$ a discrete subgroup which is large in a certain sense, e.g. of finite covolume. We are interested in studying holomorphic vector bundles over the quotient manifold $X = G/\Gamma$. In general (i.e. if $G$ is non-commutative) these manifolds are not algebraic and not even Kähler.

Thus in order to understand these manifolds it is necessary to employ methods other than those usually applied in algebraic and Kähler geometry. These methods are in particular group-theoretical methods like methode from Lie theory, representation theory, and the theory of algebraic and arithmetic groups. In order to be able to employ these group-theoretical methods we restrict our attention to holomorphic vector bundles which carry a flat holomorphic connection. Any such bundle is induced by a representation of the fundamental group. Special emphasis is put on what we call *essentially antiholomorphic* representations. For irreducible lattices $\Gamma$ in semisimple complex Lie groups $G$ of rank at least two it is an easy consequence of the celebrated arithmeticity results of Margulis et al. [7] that every flat holomorphic vector bundle over $G/\Gamma$ is isomorphic as a holomorphic vector bundle to a flat bundle induced by an essentially antiholomorphic representation. For lattices in arbitrary complex Lie groups there are always non-trivial essentially antiholomorphic representations. Moreover flat vector bundles induced by antiholomorphic representations arise naturally as higher direct images sheaves of the structure sheaf for fibrations of complex parallelizable manifolds. Therefore these results are useful to calculate Dolbeault cohomology groups via Leray spectral sequences associated to certain fibrations. In a separate paper [14] we will exploit this to deduce results on cohomology groups, deformations and the Albanese variety of complex parallellizable manifolds.

We give a precise criterion determining when a flat vector bundle induced by an essentially antiholomorphic representation is trivial (thm. 6.6.) and moreover classify these bundles completely up to isomorphism as holomorphic vector bundles in the case $G = G'$ (thm. 6.3).

We then proceed to study sections and subbundles of flat vector bundles. Essentially we prove that in general sections exist only to the extent to which the bundle

Part of the work was done during the special year on Several Complex Variables at the MSRI. Research at MSRI is supported in part by NSF grant DMS-9022140. The author also wishes to thank the Harvard University for their invitation.

is trivial and that generically every vector subbundle of a flat vector bundle on a complex parallelizable manifold is again flat.

As a by-product we can show that every compact complex manifold of dimension $n > 0$ admits a non-trivial holomorphic vector bundle of rank at most $n$ (improving our result in [11]).

The assumption that $\Gamma$ is of finite covolume is often used in a rather indirect way. For instance we prove that, given a lattice $\Gamma$ in a complex Lie group $G$, the group $G$ is always generated by those connected commutative complex Lie subgroups $A$ for which $A/(A \cap \Gamma)$ is compact (see prop. 4.1). We will also use the fact that $G/\Gamma$ carries no non-constant plurisubharmonic function.

A part of the results for the special case where $G/\Gamma$ is compact is already contained in the Habilitationsschrift of the author [12].

## 2. Generalities

Let $G$ be a connected complex Lie group and $\Gamma$ a discrete subgroup. Then $X = G/\Gamma$ is a complex parallelizable manifold, i.e. a complex manifold whose tangent bundle is holomorphically trivial. Conversely every compact complex parallelizable manifold with holomorphiclally trivial tangent bundle can be realized as a quotient of a complex Lie group by a discrete subgroup [10]. If we assume $G$ to be simply-connected (this is often convenient), then $\Gamma$ is isomorphic to the fundamental group of $X$. For any representation $\rho$ of $\Gamma$ in a complex Lie group $H$ we obtain a flat $H$-prinicpal bundle by $E = G \times H/\Gamma$ with $\Gamma$ acting on $G$ by deck transformations and on $H$ via $\rho$. Such a bundle admits a trivialization compatible with the flat connection if and only if the representation is trivial. However, we will be concerned whether there exists any trivialization of $E$ as $H$-principal bundle (not necessarily compatible with the flat connection). Such a trivialization exists if and only if there is a $\Gamma$-equivariant holomorphic map $\phi$ from $G$ to $H$ (i.e. $\phi(g\gamma) = \phi(g) \cdot \rho(\gamma)$ for all $g \in G$, $\gamma \in \Gamma$). More general, two representations $\rho, \tilde\rho : G \to H$ induce bundles which are isomorphic as holomorphic $H$-principal bundles if and only if there exists a holomorphic map $\phi : G \to H$ such that $\phi(g\gamma) = \rho(\gamma)^{-1}\phi(g)\tilde\rho(\gamma)$ for all $g \in G$ and $\gamma \in \Gamma$.

If $\Gamma$ is a "large" discrete subgroup of $G$ there are few holomorphic mappings from $X = G/\Gamma$. Especially, if $\Gamma$ is a lattice in a complex Lie group $G$ the quotient manifold $X = G/\Gamma$ has the following properties:

1. (Iwamoto's version of the Borel density theorem ([4])). Every linear representation $\rho : G \to GL(n, \mathbb{C})$ the Zariski-closures of $\rho(\Gamma)$ and $\rho(G)$ coincide. This implies that every $G$-equivariant holomorphic map from $X$ to a projective space is constant.
2. Every holomorphic and moreover every plurisubharmonic function on $X = G/\Gamma$ is constant (see e.g. claim 3.3 below).

## 3. A triviality criterion

**Proposition 1.** *Let $G$ be a connected complex Lie group, $H$ a complex Lie group and $\Gamma \subset G$ a discrete subgroup. Let $\rho : \Gamma \to H$ be a group homomorphism. Assume that every holomorphic function on $G/\Gamma$ is constant.*

*Then the $H$-principal bundle $E \to G/\Gamma$ induced by $\rho$ is holomorphically trivial if and only if $\rho$ extends to a holomorphic group homomorphism $\tilde{\rho} : G \to H$.*

*Proof.* Assume that there is a trivializing map $\phi : G \to H$, i.e. a holomorphic map $\phi : G \to H$ such that $\phi(g\gamma) = \phi(g)\rho(\gamma)$ for all $g \in G$ and $\gamma \in \Gamma$. Upon replacing $\phi$ by $\tilde{\phi}(g) \stackrel{def}{=} \phi(e)^{-1}\phi(g)$, we may assume that $\phi(e) = e$. Define $\alpha : G \times G \to H$ by

$$\alpha(g_1, g_2) = \phi(g_1)\phi(g_2)\phi(g_1 g_2)^{-1}.$$

Then $\alpha(g_1, e) = \alpha(e, g_2) = e$ and $\alpha(g_1, g_2\gamma) = \alpha(g_1, g_2)$ for all $g_1, g_2 \in G$ and $\gamma \in G$. We consider now the induced maps $\alpha_g : G/\Gamma \to H$ given by $\alpha_g : x\Gamma \mapsto \alpha(g, x)$. Since $\alpha_e \equiv e$, it is clear that all the maps $\alpha_g$ are homotopic to a constant map. Hence these maps may be lifted to the universal covering $\tilde{H}$ of $H$, i.e. there exist holomorphic maps $\tilde{\alpha}_g : G/\Gamma \to \tilde{H}$ such that $\alpha_g = \pi \circ \tilde{\alpha}_g$ where $\pi : \tilde{H} \to H$ is the universal covering map. But $\tilde{H}$ is a simply-connected complex Lie group and therefore Stein. Hence the maps $\tilde{\alpha}_g$ are constant. It follows that the maps $\alpha_g$ are likewise constant. Thus $\alpha(g, x) = \alpha(g, e) = e$ for all $g, x \in G$. □

## 4. Bounded representations

Here we study bounded representations, i.e. representations with relatively compact image. First we want to mention some ways in which bounded representations arise.

1. If $\Gamma$ is a lattice in a simply-connected complex Lie group $G$, then $\Gamma$ is finitely generated and $G$ is linear. In this case a theorem of Malcev [6] implies that $\Gamma$ is *residually finite*, i.e. for every element $\gamma \in \Gamma$ there exists a finite group $F$ and a group homomorphism $\rho : \Gamma \to F$ such that $\rho(\gamma) \neq e$. Since every finite group embeds in some linear group, this yields many representations of $\Gamma$ with finite image.
2. If $H_1(G/\Gamma) \simeq \Gamma/\Gamma'$ is non-trivial, $\Gamma$ admits many group homomorphisms to $S^1 = \{z \in \mathbb{C}^* : |z| = 1\}$.
3. Let $K$ be a number field, $\mathcal{O}_K$ its ring of algebraic integers, $S$ a semisimple $K$-group and $\mathcal{S}$ resp. $\mathcal{T}$ the set of all archimedean valuations $v$ such that $G$ is $K_v$-isotropic resp. $K_v$-anisotropic. Assume that $K_v \simeq \mathbb{C}$ for all $v \in \mathcal{S}$ and that $\mathcal{T}$ is not empty. Then $G = \Pi_{v \in \mathcal{S}} S(K_v)$ is a complex Lie group, $U = \Pi_{v \in \mathcal{T}} S(K_v)$ is a compact real Lie group and $\Gamma = S(\mathcal{O}_K)$ is a lattice in $G$ which admits an injective group homomorphism to $U$. Now every representation of $U$ induces a representation of $\Gamma$ with relatively-compact image.

**Theorem 1.** *Let $G$ be a connected complex Lie group, $R$ its radical, $\Gamma$ a discrete subgroup of $G$, $H$ a Stein Lie group and $\rho_i : \Gamma \to H$ group homomorphisms with relatively compact image.*

*Assume that $Ad(\Gamma)$ and $Ad(G)$ have the same Zariski-closure in $GL(\mathcal{L}ie\, G)$ and that $\Gamma \cap R$ is cocompact in $R$.*

*Let $E_1$ and $E_2$ denote the $H$-principal bundles over $G/\Gamma$ induced by $\rho_1$ resp. $\rho_2$. Then $E_1 \sim E_2$ if and only if $\rho_1$ and $\rho_2$ are conjugate by an element $h \in H$.*

**Remark 1.** *If $\Gamma$ is a lattice in a complex Lie group $G$ with radical $R$, then $R/(R \cap \Gamma)$ is compact ([8]) and the Zariski-closures of $Ad(\Gamma)$ and $Ad(G)$ in $GL(\mathcal{L}ie\,G)$ do coincide ([2]).*

*Proof.* Recall that $E_1 \sim E_2$ holds iff there exists a holomorphic map $\phi: G \to H$ such that

$$\phi(g\gamma) = \rho_1(\gamma)^{-1}\phi(g)\rho_2(\gamma) \tag{1}$$

for all $g \in G$, $\gamma \in \Gamma$.

Let $K_i$ denote the closure of $\rho_i(\Gamma)$ in $H$. The sets $K_i$ are compact subgroups. Since $H$ is Stein, it admits a strictly plurisubharmonic exhaustion function $\tau$. We may assume that $\tau$ is invariant under the $K_1 \times K_2$-action on $H$ given by $(k_1, k_2): h \mapsto k_1^{-1}hk_2$ (because we can replace $\tau$ by the function obtained by averaging $\tau$ over the $K_1 \times K_2$-orbits). Then (1) implies $\tau(\phi(g\gamma)) = \tau(\phi(g))$. Thus we obtain a plurisubharmonic function on $G/\Gamma$.

**Claim 1.** *Under the assumptions of the theorem every plurisubharmonic function on $G/\Gamma$ is constant.*

*Proof.* Since $R/(R \cap \Gamma)$ is compact, every plurisubharmonic function on $G/\Gamma$ is a pull-back from $(G/R)/(R\Gamma/\Gamma)$. Now $G/R$ is semisimple and $\Gamma/(\Gamma \cap R)$ is Zariski-dense in $G/R$. Hence the assertion follows from a theorem of Berteloot and Oeljeklaus [1]. $\square$

Thus $g \mapsto \tau(\phi(g))$ is constant. Since $\tau$ is strictly plurisubharmonic, this implies that $\phi$ is constant. $\square$

**Corollary 1.** *Let $G$, $\Gamma$ and $H$ be as in the above theorem and let $\rho: \Gamma \to H$ be a group homomorphism with relatively compact image.*

*Then the induced $H$-principal bundle over $G/\Gamma$ is holomorphically trivial if and only if $\rho \equiv e$.*

## 5. Subgroups with a bounded orbit

Given a Lie group $G$ and a lattice $\Gamma$ in $G$ we will show that there are many Lie subgroups $H$ of $G$ such that the $H$-orbit through $e\Gamma$ is relatively compact in $G/\Gamma$.

**Proposition 2.** *Let $G$ be a connected complex Lie group, $\Gamma$ a lattice and $G_0$ the subgroup of $G$ generated by all connected commutative complex Lie subgroups $A$ of $G$ with $A/(A \cap \Gamma)$ compact.*
*Then $G_0 = G$.*

*Proof.* If $G$ is commutative, then $G/\Gamma$ is a topological group with finite volume and therefore compact. In the general case we argue by induction over $\dim(G)$. Thus let $G$, $\Gamma$ and $G_0$ be as above and assume that the proposition is valid for all Lie groups of lower dimension. Let $\gamma \in \Gamma \setminus Z$ (where $Z$ denotes the center of $G$) and let $C^0(\gamma)$ denote the connected component of the centralizer of $\gamma$ in $G$. Then $\dim C^0(\gamma) < \dim G$ and $C^0(\gamma) \cap \Gamma$ is a lattice in $C^0(\gamma)$ ([9],Lemma 1.14). Hence $C^0(\gamma) \subset G_0$ by the induction hypothesis. This is equivalent to the assertion that $V^{Ad(\gamma)} \subset \mathcal{L}ie(G_0)$ where $V = \mathcal{L}ie\,G$. Let $H \subset GL(V)$ denote the Zariski-closure of $Ad(G)$. Then $H$ is also the Zariski-closure of $Ad(\Gamma)$. Let $k$ denote the generic dimension of $V^h = \{v : hv = v\}$ for $h \in H$ and $\Omega$ the set of all $h \in H$ for which $\dim V^h = k$. Then $\Omega$ is a Zariski-open subset of $H$. Now $V^h \subset \mathcal{L}ie(G_0)$ for all $h \in \Gamma \cap \Omega$ and therefore for all $h \in \Omega$. It follows that $v \in \mathcal{L}ie\, G_0$ for all $v \in Ad^{-1}(\Omega)$. Since $Ad^{-1}(\Omega)$ is a non-empty open subset of $\mathcal{L}ie\, G$, this implies that $G = G_0$. $\square$

We will now deduce another variant of this theme, this time strengthening the assumption on the subgroups (requiring unipotency) while relaxing the assumption on the orbits (only relatively-compact instead of compact).

**Proposition 3.** *Let $G$ be a simply-connected complex Lie group and $\Gamma$ a lattice.*

*Let $G_1$ denote the subgroup of $G$ generated by all unipotent subgroups $U \subset G'$ for which there exists a compact subset $F \subset G$ such that $U \subset F\Gamma'$.*

*Then $G_1 = G'$.*

**Remark 2.** *For a simply-connected complex Lie group $G$ the commutator group $G'$ carries a unique structure of a complex linear-algebraic group. Hence it makes sense to speak about unipotent subgroups of $G'$.*

*Proof.* First we note that we may assume $R \cap G' = \{e\}$. because $R \cap G'$ is a normal unipotent subgroup with compact orbits. This assumption implies that $R$ is central and $G = S \times R$ with $S$ semisimple. Let $\pi : G \to G/R \simeq S$ be the natural projection.

**Claim 2.** *Let $H$ be a one-dimensional unipotent subgroup of $S$ such that $H \cap \pi(\Gamma) \neq \{e\}$.*

*Then $H \subset G_1$.*

*Proof.* Let $u \in S$, $r \in R$ such that $ur \in \Gamma$ and $\pi(u) \in H \setminus \{e\}$. Let $C$ denote the connected component of the centralizer of $u$ in $S$. Clearly, $H$ is contained in the center of $C$. Now $C \cdot R$ is the centralizer of $ur$ in $G$. Hence the $CR$-orbit through $e\Gamma$ is closed. Let $N$ denote the nilradical of $C$. Then $NR$ is the nilradical of $CR$. It follows that the $NR$-orbit through $e\Gamma$ is compact. Now $H$ is central in $C$, hence $H \subset N$ and consequently $H \subset G_1$. $\square$

Next let $V$ denote the subgroup of $G/R$ generated by all one-dimensional unipotent subgroups having non-trivial intersection with $\pi(\Gamma)$. Then $V$ is normalized by $\pi(\Gamma)$. The Borel Density theorem thus implies that $V$ is a normal subgroup of $S$. It follows that $S = S_0 \times V$ for some semisimple complex Lie subgroup $S_0 \subset S$.

**Claim 3.** *The $S_0$-orbit through $e\Gamma$ in $G/\Gamma$ is relatively compact.*

*Proof.* By the criterion of Kazdan-Margulis (see [9], Cor. 11.12) we have to show that given a sequence $s_n \in S_0$ there is no sequence $\gamma_n \in \Gamma \setminus \{e\}$ such that $\lim s_n \gamma_n s_n^{-1} = e$. Assume there are such sequences. Let $u_n \in S$ and $r_n \in R$ such that $u_n r_n = \gamma_n$. Then $\lim s_n u_n s_n^{-1} = e$ (because $R$ is central). Since $\pi(\Gamma)$ is a lattice in $S$, it follows that $u_n$ is unipotent for $n$ sufficiently large ([9], Cor. 11.18). But in this case $u_n \in V$ by construction and this implies that $s_n u_n s_n^{-1} = u_n$ for all $n$. Contradiction! $\square$

Finally note that $S_0$ is a semisimple complex Lie group and therefore generated by its unipotent subgroups. $\square$

## 6. Antiholomorphic maps and actions of unipotent groups

Here we prove an auxiliary result needed later on.

**Proposition 4.** *Let $Z$ be a quasi-affine variety, $\bar{H}$ a connected commutative linear-algebraic group (both defined over $\mathbb{C}$) and $H$ a connected Zariski-dense complex Lie subgroup of $\bar{H}$. Assume that there is a regular action $\mu : \bar{H} \times Z \to Z$, a map $\alpha : H \to Z$ with relatively compact image and an antiholomorphic map $\phi : H \to Z$ such that*

$$\alpha(h) = \mu(h)(\phi(h)) \tag{2}$$

*for all $h \in H$. Let $U$ denote the unipotent radical of $\bar{H}$.*

*Then both $\alpha(H)$ and $\phi(H)$ are contained in a single $\mu(\bar{H})$-orbit $W$, $U$ acts trivially on this orbit and the maps $\alpha$, $\phi$ are homomorphisms of real Lie groups from $H$ to $W$ with respect to the natural group structure on $W$ with $\phi(e)$ as neutral element.*

*Proof.* We start with a discussion of *semi-invariant* functions on $Z$.

**Definition 1.** *Let $Z$ be a variety, $G$ an algebraic group acting on $Z$. A regular function $f$ on $Z$ is called* semi-invariant *if there exists a character $\chi$ of $G$ such that $f(gz) = \chi(g)f(z)$ for all $z \in Z$ and $g \in G$.*

Now let $f$ be a semi-invariant for the $\bar{H}$-action on $Z$. Then

$$f(\alpha(h)) = \chi(h) f(\phi(h)) \tag{3}$$

for a character $\chi$ of $\bar{H}$ and all $h \in H$. This implies

$$|f(\alpha(h))|^2 = |\chi(h) f(\phi(h))|^2 = |\chi(h) \overline{f(\phi(h))}|^2$$

The left side is bounded while the right side is the absolute value of a holomorphic function. Hence there is a constant $c$ such that $|f \circ \alpha| \equiv c$ and $\chi \cdot \overline{f \circ \phi} \equiv c$. It follows that $f \circ \alpha$ and $f \circ \phi$ vanish either everywhere or nowhere.

**Lemma 1.** *Let $W$ be a quasi-affine algebraic variety and $G$ a connected solvable linear-algebraic group acting regularly on $W$. Assume that for every semi-invariant $f$ the zero-set $V(f)$ is either empty or the whole of $W$.*

*Then $G$ acts transitively on $W$.*

*Proof.* Assume the contrary. Then $W$ must contain a proper invariant algebraic subvariety $Y$. Now the ideal $I_Y$ is a non-trivial invariant subvectorspace of the space of regular functions $\mathbb{C}[W]$. Recall that every $f \in \mathbb{C}[W]$ is contained in a finite-dimensional subvectorspace. Using this fact the theorem of Lie implies that $I_Y$ contains a one-dimensional invariant subvectorspace $S$. Now any $f \in S \setminus \{0\}$ is a semi-invariant on $W$ vanishing on $Y$ but not vanishing everywhere. Contradiction! □

Applying this lemma to our situation we may conclude that both $\alpha(H)$ and $\phi(H)$ are contained in a single $\bar{H}$-orbit $W$.

Since $\bar{H}$ is commutative, the homogeneous $\bar{H}$-space $W$ has a canonical structure as a commutative group which is unique up to the choice of the neutral element.

**Claim 4.** *The $U$-action on $W$ is trivial.*

*Proof.* If not, there is a non-trivial regular $U$-equivariant map $\tau : W \to \mathbb{C}$. Using the assumption that $H$ is Zariski-dense in $\bar{H}$, we may find a one-parameter subgroup $\gamma(t)$ in $H$ such that $\tau(\mu(\gamma(t))(x)) = t + \tau(x)$ for $x \in W$. But then $t \mapsto \tau(\alpha(\gamma(t)))$ is a bounded function on $\mathbb{C}$ which may be represented as

$$\tau(\alpha(\gamma(t))) = \tau(\mu(\gamma(t))(\phi(t))) = t + \underbrace{\tau(\phi(\gamma(t)))}_{antiholo.}. \tag{4}$$

This is a contradiction, because the left side is bounded while the right side is a non-constant harmonic function on $\mathbb{C}$. Thus the $U$-action must have been trivial. $\square$

Now $W$ is a homogeneous space of the reductive commutative group $\bar{H}/U$. By choosing a point in $W$ as neutral element, $W$ inherits a structure as reductive commutative algebraic group. Let us choose $\phi(e) = \alpha(e)$ as neutral element. Then every character of $W$ is an semi-invariant for the $\bar{H}$-action and consequently our previous considerations imply that $\chi \circ \phi$ is a real Lie group homomorphism from $H$ to $\mathbb{C}^*$ for every character $\chi$ of $W$. It follows that $\phi$ (and hence also $\alpha$) are real Lie group homomorphisms. $\square$

## 7. Essentially antiholomorphic representations

**Definition 2.** *Let $G$ and $H$ be complex Lie groups, and $\Gamma \subset G$ a discrete subgroup. A group homomorphism $\rho : \Gamma \to H$ is called* essentially antiholomorphic *if there exists an antiholomorphic Lie group homomorphism $\zeta : G \to H$ and a map $\xi : \Gamma \to H$ with relatively compact image $\xi(\Gamma) \subset H$ such that*

$$\rho(\gamma) = \zeta(\gamma) \cdot \xi(\gamma) \tag{5}$$

*for all $\gamma \in \Gamma$.*

**Proposition 5.** *Let $G$ be a connected complex Lie group, $\Gamma$ a lattice, $H$ a Stein Lie group and $\rho : \Gamma \to H$ an essentially antiholomorphic representation with maps $\zeta$, $\xi$ given as above.*

*Then*

1. *$\xi : \Gamma \to H$ is a group homomorphism.*
2. *Both $\xi$ and $\zeta$ are uniquely determined by $\rho$.*
3. *$\xi(\gamma)$ and $\zeta(g)$ commute for all $g \in G$, $\gamma \in \Gamma$.*

*Proof.* The equation (5) combined with $\rho(\gamma\delta) = \rho(\gamma)\rho(\delta)$ implies

$$\zeta(\delta)^{-1}\xi(\gamma)\zeta(\delta) = \xi(\gamma\delta)\xi(\delta)^{-1} \tag{6}$$

for all $\gamma, \delta \in \Gamma$.

Let $A$ be a connected complex Lie subgroup of $G$ for which there exists a compact set $F \subset G$ such that $A \subset F\Gamma$. For $\gamma \in \Gamma$ we define an antiholomorphic map $\phi_\gamma : A \to H$ by

$$\phi_\gamma(a) = \zeta(a)\xi(\gamma)\zeta(a)^{-1}.$$

If $\overline{\xi(\Gamma)} = K$, then $\phi_\gamma(A) \subset \zeta(F) \cdot K \cdot K^{-1} \cdot \zeta(F)$. Since the latter set is compact, it follows that $\phi_\gamma$ is constant. This implies that $\zeta(A)$ commutes with $\xi(\Gamma)$. Since $G$ is generated as a group by such subgroups $A$ (prop. 2), this yields the third assertion. The first assertion is an immediate corollary. To check unicity, let $\zeta$, $\tilde\zeta$ be antiholomorphic representations and $\xi$, $\tilde\xi$ be bounded maps such that $\zeta\xi = \tilde\zeta\tilde\xi$. This implies $\tilde\zeta^{-1}\zeta = \tilde\xi\xi^{-1}$. This is now a bounded antiholomorphic map, hence constant, and in fact constant with value $e$. Therefore $\tilde\zeta = \zeta$ and $\tilde\xi = \xi$. □

**Theorem 2.** *Let $G$ be a connected complex Lie group, $H$ a linear complex Lie group, $\Gamma$ a lattice and $\rho, \tilde\rho : \Gamma \to H$ essentially antiholomorphic representations. Assume that the induced flat bundles are isomorphic as holomorphic $H$-bundles.*

*Then there exists a constant $c \in H$ such that $\rho|_{\Gamma \cap G'} = c\,(\tilde\rho|_{\Gamma \cap G'})\,c^{-1}$.*

*Proof.* There is no loss in generality in assuming that $G$ is simply-connected. The assumption of the induced flat bundles being holomorphically equivalent is equivalent to the statement that there exists a holomorphic map $\phi : G \to H$ such that

$$\phi(g\gamma) = \xi(\gamma)^{-1}\zeta(\gamma)^{-1}\phi(g)\tilde\zeta(\gamma)\tilde\xi(\gamma) \tag{7}$$

for all $g \in G$ and $\gamma \in \Gamma$ where $\rho = \zeta \cdot \xi$ resp. $\tilde\rho = \tilde\zeta \cdot \tilde\xi$ are the respective decompositions as products of an antiholomorphic and a bounded factor. We define a map $\alpha : G \to H$ by

$$\alpha(g) = \zeta(g)\phi(g)\tilde\zeta(g)^{-1} \tag{8}$$

Then

$$\begin{aligned}\alpha(g\gamma) &= \zeta(g)\zeta(\gamma)\xi(\gamma)^{-1}\zeta(\gamma)^{-1}\phi(g)\tilde\zeta(\gamma)\tilde\xi(\gamma)\tilde\zeta(\gamma)^{-1}\tilde\zeta(g)^{-1} \\ &= \xi(\gamma)^{-1}\zeta(g)\phi(g)\tilde\zeta(g)^{-1}\tilde\xi(\gamma) = \xi(\gamma)^{-1}\alpha(g)\tilde\xi(\gamma)\end{aligned} \tag{9}$$

(Recall that $\zeta(G)$ and $\xi(\Gamma)$ commute).

**Claim 5.** *Let $A$ be a connected commutative complex Lie subgroup of $G$ such that there exists a compact subset $F_A \subset G$ such that $A \subset F_A\Gamma$. Let $c = \phi(e) = \alpha(e) \in H$ and $W = \mu(G)(c)$ where $\mu$ is the action given by*

$$\mu(g) : x \mapsto \zeta(g)^{-1}x\tilde\zeta(g) \tag{10}$$

*Let $\eta$ denote the $\Gamma$-action on $H$ given by*

$$\eta(\gamma) : x \mapsto \xi(\gamma)^{-1}x\tilde\xi(\gamma) \tag{11}$$

*Then $\alpha(A) \subset W$, $\phi(A) \subset W$ and $\eta(A \cap \Gamma)$ stabilizes $W$. If in addition $A$ is contained in $G'$, then $\mu(A)$ acts trivially on $W$.*

*Proof.* Fix an antiholomorphic embedding $\sigma : H \hookrightarrow GL(n, \mathbb{C})$ and let $\hat\mu$ denote the $G$-action on $GL(n, \mathbb{C})$ given by

$$\hat\mu(g) : x \mapsto \sigma \circ \zeta(g^{-1}) \cdot x \cdot \sigma \circ \tilde\zeta(g).$$

Thus $\hat\mu$ embedds $G$ holomorphically in $GL(n, \mathbb{C}) \times GL(n, \mathbb{C})$ acting on $GL(n, \mathbb{C})$ by left and right multiplication. Now an application of prop. 4 (with $Z = GL(n, \mathbb{C})$) and the Zariski-closure of $\hat{\mu}(H)$ in $GL(n,\mathbb{C}) \times GL(n,\mathbb{C})$ as $\bar{A}$) yields $\alpha(A) \subset W$ and $\phi(A) \subset W$.

By definition of $\alpha$ the inclusion $(A \cap \Gamma) \subset W$ implies that $\eta(\gamma)(c) \in W$ for $\gamma \in A \cap \Gamma$. Since $\eta(\Gamma)$ commutes with $\mu(G)$, it follows that $\eta(A \cap \Gamma)$ stabilizes $W$. Finally, since $G$ was assumed to be simply-connected, this implies that $G'$ is unipotent. Hence the last assertion also follows from prop. 4. □

**Claim 6.** *The group $\mu(G')$ acts trivially on $W$ and in particular stabilizes the point $\phi(e)$.*

This is immediate, because $G'$ is generated by unipotent subgroups with a bounded orbit (prop. 3).

Next we discuss 'translates' of $\phi$. For $x \in G$ let $\phi_x$ denote the map defined by $\phi_x(g) = \phi(xg)$. Then $\phi_x : G \to H$ is a holomorphic map fulfilling (7). Thus we may apply the above considerations and conclude that $\mu(G')(\phi_x(e)) = \phi_x(e)$. Since $\phi_x(e) = \phi(x)$, it follows that $\mu(G')$ stabilizes $\phi(G)$ pointwise. Hence

$$\phi(g\gamma) = \eta(\gamma)(\phi(g)) \tag{12}$$

for $g \in G$, $\gamma \in G' \cap \Gamma$. By assumption $\xi(\Gamma)$ and $\tilde{\xi}(\Gamma)$ are contained in compact subgroups $K$ resp. $\tilde{K}$ of $H$. Since $H$ is Stein, it admits a strictly plurisubharmonic exhaustion function $\tau$. By averaging with respect to Haar measure we may assume that $\tau$ is invariant with respect to the given $K \times \tilde{K}$-action on $H$. Then we obtain

$$\tau(\phi(g\gamma)) = \tau(\phi(g)) \tag{13}$$

for $g \in G$, $\gamma \in G' \cap \Gamma$. Now every plurisubharmonic function on $G'/(G' \cap \Gamma)$ is constant (see [14]). Hence $\phi$ is constant along the $G'$-orbits on $G$. In particular $\phi|_{G'}$ is constant and this yields the assertion with $c = \phi(e)$. □

**Theorem 3.** *Let $G$ be a connected complex Lie group, $\Gamma \subset G$ a lattice, $H$ a complex linear-algebraic group and $\rho : \Gamma \to H$ an essentially antiholomorphic group homomorphism.*

*Then the flat bundle over $G/\Gamma$ defined by $\rho$ is trivial if and only if the following conditions are fulfilled:*
 1. *The Zariski-closure $A$ of $\rho(\Gamma)$ in $H$ is isomorphic to $(\mathbb{C}^*)^g$ for some $g \geq 0$.*
 2. *Let $\rho = \zeta \cdot \xi$ denote the decxomposition of $\rho$ into the antiholomorphic part $\zeta$ and the bounded part $\xi$. Let $K$ denote the maximal compact subgroup of $A$ (which is a real form of $A$) and let $\tau$ denote complex conjugation on $A$ with respect to the real form $K$ (i.e. $\tau(w_1, \ldots, w_g) = (\bar{w}_1^{-1}, \ldots, \bar{w}_g^{-1})$).*

*Then $\rho(\gamma) \cdot \tau(\zeta(\gamma)) = e$ for all $\gamma \in \Gamma$.*

*Proof.* By prop. 1 triviality of the induced bundle is equivalent to the existence of a holomorphic Lie group homomorphism $\phi : G \to H$ such that $\phi|_\Gamma = \rho$. Hence one direction of the statement is clear: If the criterion is fulfilled, then $\phi : g \mapsto \tau(\zeta(g)^{-1})$ is a holomorphic Lie group homomorphism from $G$ to $A \subset H$ such that $\phi|_\Gamma = \rho$.

Let us assume that the bundle is trivial and that $\phi : G \to H$ is the trivializing holomorphic Lie group homomorphism. Then $\phi|_{G'}$ is constant by the preceding theorem. This implies that $\rho(\gamma) = e$ for all $\gamma \in G' \cap \Gamma$. Hence $\rho(\Gamma)$ is commutative.

Let $A$ denote the Zariski-closure of $\rho(\Gamma)$ in $H$. Then $H/A$ is a quasi-affine variety (because $A$ is commutative). Thus $\phi : G \to H$ induces a holomorphic map from $G/\Gamma$ to $H/A$ which must be constant. We may assume that $\phi(e) = e$. Then $\phi(G) \subset A$. Hence we may assume that $H = A$, i.e. we may assume that $H$ is a commutative linear-algebraic group. Thus $H \simeq (\mathbb{C})^d \times (\mathbb{C}^*)^g$ and we may discuss these factors separately.

Let us first discuss the case $H = \mathbb{C}$. There does not exist any relatively compact subgroup in $(\mathbb{C}, +)$. Therefore, if $\rho = \zeta \cdot \xi$ is the decomposition of $\rho$ into the antiholomorphic part $\zeta$ and the bounded part $\xi$, then $\xi \equiv 0$, i.e. $\rho = \zeta$. It follows that $\phi - \rho$ induces a pluriharmonic function on $G/\Gamma$. But every pluriharmonic function on $G/\Gamma$ is constant. This forces both $\phi$ and $\rho$ to be constant, i.e. $\rho \equiv e$.

Thus we have seen that the Zariski-closure $A$ of $\rho(\Gamma)$ cannot contain a factor isomorphic to $\mathbb{C}$, i.e. $A$ must be isomorphic to some $(\mathbb{C}^*)^g$. Let $K \simeq (S^1)^g$ denote the maximal compaxt subgroup of $A$. This is a real form, hence there is a map $\tau$ given by complex conjugation on $A$ with respect to $K$. (In coordinates this is just $(\bar{w}_1^{-1}, \ldots, \bar{w}_g^{-1})$.) Consider now the map $\alpha : g \mapsto \phi(g)\tau(\zeta(g))$. Note that

$$\alpha(\gamma) = \underbrace{\zeta(\gamma)\tau(\zeta(\gamma))}_{\in K} \underbrace{\xi(\gamma)}_{\in K}.$$

Furthermore, $\Gamma$ being a lattice in the Lie group $G$ implies that the closure (in the euclidean topology) of $G'\Gamma$ in $G$ is cocompact. This implies that there exists a compact subset $C \subset G$ such that $C \cdot G' \cdot \Gamma$ is dense in $G$. By the above considerations $\alpha$ maps this dense subset of $G$ into a compact subset of $A$, namely $\alpha(C) \cdot K$. Since $\alpha$ is a group homomorphism, it follows that $\alpha(G) \subset K$. But $\alpha$ is also a holomorphic map. Hence $\alpha(G) \subset K$ implies that $\alpha$ is constant, i.e. $\alpha \equiv e$. Thus $\alpha(\gamma) = \rho(\gamma) \cdot \tau(\zeta(\gamma)) = e$ for all $\gamma \in \Gamma$. □

In the special case where $\rho$ is bounded (i.e. $\zeta \equiv e$) this reproves cor.1 of thm. 3.1.

In an other special case, $\rho = \zeta|_\Gamma$ and $\xi \equiv e$, we obtain the following criterion.

**Corollary 2.** *Let $G$ be a complex Lie group, $\Gamma \subset G$ a lattice, $H$ a complex linear-algebraic group and $\rho : G \to H$ an antiholomorphic Lie group homomorphism.*

*Then the flat bundle over $G/\Gamma$ defined by $\rho|_\Gamma$ is trivial if and only if there exists a number $n$ and holomorphic Lie group homomorphisms $\rho_1 : G \to (\mathbb{C}^*)^n$ and $\rho_2 : (\mathbb{C}^*)^n \to H$ such that $\rho = \rho_2 \circ \zeta \circ \rho_1$ where $\zeta$ denotes complex conjugation on $T = (\mathbb{C}^*)^n$ and $\rho_1(\Gamma) \subset (\mathbb{R}^*)^n$.*

## 8. Vector Bundles

To transfer our results on principal bundles to vector bundles, we need the following generalized Schur lemma.

**Lemma 2.** *Let $\Gamma$ be a group, $V$ a vector space and $\zeta, \xi : \Gamma \to GL(V)$ group homomorphisms such that $\zeta(\gamma)$, $\xi(\delta)$ commute for all $\gamma, \delta \in \Gamma$. Assume that $V$ is an irreducible $\Gamma$-module for $\rho = \zeta \cdot \xi$.*

*Then there exists vector spaces $V_1$, $V_2$ and irreducible representations $\zeta_0 : \Gamma \to GL(V_1)$, $\xi_0 : \Gamma \to GL(V_2)$ such that $(V, \rho) \simeq (V_1, \zeta_0) \otimes (V_2, \xi_0)$.*

**Lemma 3.** *Let $A \subset GL(n, \mathbb{C})$ be a connected commutative complex Lie subgroup, $\Lambda \subset A$ a lattice, $Q$ a projective manifold with $b_1(Q) = 0$ on which $A$ acts holomorphically and $x \in Q$. Assume that for every $\lambda \in \Lambda$ there exists a sequence of natural numbers $n_k$ such that $\lim n_k = \infty$ and $\lim(\lambda^{n_k}(x)) = x$.*

*Then $x$ is a fixed point for the $A$-action on $Q$.*

*Proof.* The assumptions on $Q$ imply that $Aut(Q)^0$ is a linear-algebraic group. Let $H$ denote the Zariski-closure of $A$ in $Aut(Q)^0$. Assume that $x$ is not a fixed point for the $A$-action. The $H$-orbit through $x$ is a locally closed subset of $Q$ and isomorphic to the quotient $H/I$ where $I = \{h \in h : h(x) = x\}$. This quotient $H/I$ is again a linear-algebraic group (because $H$ is commutative). Recall that on a complex Lie group $A$ every bounded holomorphic function is constant. It follows that the image of $A$ in $H/I$ under the natural morphism $\tau : A \to H/I$ can not be relatively compact. Hence there is an element $\lambda \in \Lambda$ such that the sequence $\tau(\lambda^n)$ in $H/I$ contains no convergent subsequence. Since the $H$-orbits in $Q$ are locally closed, this implies that no subsequence of $\lambda^n(x)$ can converge to $x$. This contradicts the assumptions of the lemma. Hence $x$ must be an $A$-fixed point. □

**Proposition 6.** *Let $G$ be a connected complex Lie group, $\Gamma$ a lattice, $V$ a complex vector space and $\rho : \Gamma \to GL(V)$ be an essentially antiholomorphic representation, with antiholomorphic part $\zeta$ and bounded part $\xi$. Assume that $V$ contains a $\rho(\Gamma)$-invariant subvectorspace $W$.*

*Then $W$ is already invariant under $\zeta(G)$ and $\xi(\Gamma)$.*

*Proof.* Let $k = \dim W$. We consider the induced actions on the Grassmann manifold $Q$ of $k$-dimensional subvectorspaces of $V$. For simplicity, they are also denoted by $\rho$, $\zeta$ and $\xi$. Let $x = [W] \in Q$. Since $\xi(\Gamma)$ is bounded, it is clear that for every element $\gamma \in \Gamma$ there exists a sequence of natural numbers $n_k$ such that $\lim n_k = \infty$ and $\lim \xi(\gamma^{n_k}) = e$. Then $\rho = \zeta \cdot \xi$ implies that $\lim \zeta(\gamma^{-n_k})(x) = x$. By lemma 3 it follows that $x$ is a fixed point for every connected commutative complex Lie subgroup $A \subset G$ with $A/(A \cap \Gamma)$ compact. By prop. 2 this implies that $x$ is a fixed point for $\zeta(G)$ (and hence for $\xi(\Gamma)$, too). □

**Corollary 3.** *Let $G$ be a connected complex Lie group, $\Gamma$ a lattice, $V$ a complex vector space, $\rho : \Gamma \to GL(V)$ an essentially antiholomorphic representation and $W$ a $\rho(\Gamma)$-stable subvectorspace of $V$.*

*Then the restricted representation $\rho' : \Gamma \to GL(W)$ is likewise essentially antiholomorphic.*

**Lemma 4.** *Let $\Gamma$ be a group, $\rho_i : \Gamma \to GL(V_i)$ representations on complex vector spaces for $i = 1, 1$ and assume that*

1. *Both $V_i$ are irreducible $\Gamma$-modules with respect to the representations $\rho_i$.*
2. *The Zariski-closure $H$ of $\rho_1(\Gamma)$ in $GL(V_1)$ is connected.*
3. *The image $\rho_2(\Gamma)$ is finite.*

*Then $V_1 \otimes V_2$ is an irreducible $\Gamma$-module with respect to $\rho_1 \otimes \rho_2$.*

*Proof.* Let $\Gamma_0 = \ker \rho_2$. Since $\Gamma/\Gamma_0$ is finite and $H$ is connected, it is clear that $\rho_1(\Gamma_0)$ is Zariski-dense in $GL(V_1)$. Hence $V_1$ is an irreducible $\Gamma_0$-module while $\Gamma_0$ acts trivially on $V_2$. It follows that every $\Gamma_0$-invariant subvectorspace of $V_1 \otimes V_2$ has the form $V_1 \otimes W$ for some subvectorspace $W \subset V_2$. Clearly, such a $V_1 \otimes W$ is $\Gamma$-invariant only if $W = \{0\}$ or $W = V_2$. Thus $V_1 \otimes V_2$ is an irreducible $\Gamma$-module.　□

**Corollary 4.** *Let $G$ be a connected complex Lie group, $\Gamma \subset G$ a lattice, $\rho_1 : G \to GL(V_1)$ an antiholomorphic representation, $\rho_2 : \Gamma \to GL(V_2)$ a representation with relatively compact and Zariski-connected image and $\rho_3 : \Gamma \to GL(V_3)$ be a representation with finite image.*

*Assume that all the representations $\rho_i$ are irreducible.*

*Then $(\rho_1|_\Gamma) \otimes \rho_2 \otimes \rho_3 : \Gamma \to GL(V_1 \otimes V_2 \otimes V_3)$ is likewise irreducible.*

## 9. The classification

We make use of Margulis' superrigidity theorem in the following form:

**Theorem 4.** *Let $S$ be a simply-connected semisimple complex Lie group and $\Gamma \subset S$ a lattice. Assume that there does not exist a normal Lie subgroup $S_0 \subset S$ such that $S_0 \simeq SL(2, \mathbb{C})$ and $S_0/(S_0 \cap \Gamma)$ is of finite volume.*

*Then there exists a compact real semisimple Lie group $K$ and a group homomorphism $j : \Gamma \to K$ such that for every complex-algebraic group $H$ and every group homomorphism $\alpha : \Gamma \to H$ there exist continuous group homomorphisms $\zeta : S \to H$, $\xi : K \to H$ and $\nu : \Gamma \to H$ such that*

1. *$\alpha(\gamma) = \zeta(\gamma) \cdot \xi(j(\gamma)) \cdot \nu(\gamma)$ for all $\gamma \in \Gamma$.*
2. *The image $\nu(\Gamma)$ is finite.*
3. *$\zeta(s)$, $\xi(k)$ and $\nu(\gamma)$ commute for every $s \in S$, $k \in K$ and $\gamma \in \Gamma$.*

Note that for every continuous group homomorphism $\zeta$ from a complex semisimple Lie group $S$ to a complex algebraic group $H$ there exists a holomorphic group homomorphism $\zeta_0$ and an antiholomorphic group homomorphism $\zeta_1$ such that $\zeta = \zeta_0 \cdot \zeta_1$.

Using Margulis' theorem and our previous results we obtain the following classification.

**Theorem 5.** *Let $S$ be a simply-connected semisimple complex Lie group and $\Gamma \subset S$ a lattice. Assume that there does not exist a normal Lie subgroup $S_0 \subset S$ such that $S_0 \simeq SL(2, \mathbb{C})$ and $S_0/(S_0 \cap \Gamma)$ is of finite volume.*

*Then there exists a compact real semisimple Lie group $K$ and a group homomorphism $j : \Gamma \to K$ such that there is a one-to-one correspondance between irreducible holomorphic vector bundles on $S/\Gamma$ admitting a flat connection and triples $(\rho, \zeta, \nu)$ where*

1. *All $\rho$, $\zeta$ and $\nu$ are irreducible representations of $\Gamma$.*
2. *$\rho$ extends to an antiholomorphic representation of $G$.*
3. *$\zeta$ fibers through a representation of $K$.*
4. *$\nu$ is a representation with finite image.*

## 10. Sections of Homogeneous Vector Bundles

We prove that a homogeneous vector bundle over a quotient $X = G/\Gamma$ of a connected complex Lie group $G$ by a lattice $\Gamma$ admits global sections only inasmuch as it is trivial.

Since we do not want to assume $G/\Gamma$ to be compact, we first have to show that the space of sections is finite-dimensional.

**Proposition 7.** *Let $G$ be a connected complex Lie group, $\Gamma$ a lattice and $E \to X = G/\Gamma$ a homogeneous vector bundle.*
*Then $\Gamma(X, E)$ is finite-dimensional.*

*Proof.* We first need a result on the structure of $G/\Gamma$.

**Claim 7.** *Let $G$ be a connected complex Lie group and $\Gamma$ a lattice. Then there exists a $G$-equivariant holomorphic surjective map $\pi$ from $X$ onto a compact complex parallelizable manifold $Y$ such that the algebraic dimension of the fibers equals zero.*

*Proof.* The algebraic reduction of $X$ maps $X$ onto a compact complex torus $T$. If the fibers have algebraic dimension larger than zero, we continue by considering the algebraic reduction of the fiber. This yields us a holomorphic surjection of $X$ onto a compact complex parallelizable manifold (which is a torus bundle over a torus). We may continue in this way and for dimension reasons we will arrive at a holomorphic surjective map from $X$ onto a compact complex parallelizable manifold for which the fibers have algebraic dimension zero. □

We will now discuss $E$ restricted to a fiber of $\pi$.

**Claim 8.** *Let $F$ be a complex manifold of algebraic dimension zero (i.e. every meromorphic function on $F$ is constant) and $E \to F$ a vector bundle of rank $r$. Then $\dim(\Gamma(F, E)) \leq r$.*

*Proof.* Assume the contrary. Then there exists a number $d$ with $1 \leq d \leq r$ and sections $\sigma_0, \ldots, \sigma_d$ such that

1. The sections $(\sigma_i)_{0 \leq i \leq d}$ are linearly independent as elements in $\Gamma(F, E)$.
2. For a generic point $x \in F$ the subvectorspace of $E_x$ spanned by $(\sigma_i(x))_{1 \leq i \leq d}$ has dimension $d$.
3. For every point $x \in F$ the subvectorspace of $E_x$ spanned by $(\sigma_i(x))_{0 \leq i \leq d}$ has dimension at most $d$.

But in this case one of the meromorphic functions
$$f_i = (\sigma_0 \wedge \sigma_1 \wedge \ldots \wedge \widehat{\sigma_i} \wedge \ldots \wedge \sigma_d)/(\sigma_1 \wedge \ldots \wedge \sigma_d)$$
must be non-constant, contradicting the assumption of $F$ having algebraic dimension zero. □

Now let $X = G/\Gamma$, $E \to X$ and $\pi : X \to Y$ as above. Since $\dim(F, E) \leq rank(E)$ for every fiber $F$ of $\pi$, we may conclude that the direct image sheaf $\pi_* E$ is finitely generated as a $\mathcal{O}_Y$-module sheaf. For homogeneity reasons it is locally free and therefore coherent. Since $Y$ is compact, it follows that $\Gamma(Y, \pi_* E) \simeq \Gamma(X, E)$ is finite-dimensional. □

Now we are in a position to prove the following structure theorem on sections of homogeneous vector bundles.

**Proposition 8.** *Let $G$ be a connected complex Lie group, $\Gamma$ a lattice and $E \to X = G/\Gamma$ a flat vector bundle.*

*Then $E$ contains a vector subbundle $E_0$ which is parallel with respect to the flat connection and trivial as a holomorphic vector bundle such that $\Gamma(X, E) = \Gamma(X, E_0)$.*

*Proof.* The sections of $E$ generate a coherent subsheaf $E_0$ of $E$. By homogeneity this subsheaf is locally free. Evidently it is invariant under the $G$-action, hence parallel with respect to the connection. There is an equivariant map from $X$ to the projective space $\mathbb{P}(\Gamma(X, E_0)^*)$ (which is finite-dimensional due to the preceding proposition). This map is constant and hence $E_0$ is trivial as a holomorphic vector bundle. □

Together with prop. 3 this enables us to completely determine the global sections in a flat vector bundle on $G/\Gamma$ which is induced by an antiholomorphic representation of $G$. We will apply this in [14] in order to calculate the dimension of Dolbeault cohomology groups by use of Leray spectral sequences.

**Corollary 5.** *Let $G$ be a connected complex Lie group, $\Gamma$ a lattice and $\rho$ be a holomorphic representation of $G$ on a complex vector space $V$. Let $V_1 = V^{G'}$ denote the vector subspace of all $v \in V$ which are fixed by $\rho(G')$. Let $\Sigma$ denote the subset of all $v \in V_1$ such that $v$ is an eigenvector with real eigenvalue for every $\rho(\gamma)$ ($\gamma \in \Gamma$) and let $V_0$ denote the subvectorspace of $V$ spanned by $\Sigma$.*

*Let $E$ and $E_0$ denote the flat vector bundles on $X = G/\Gamma$ which is induced by the representation $\bar\rho|_\Gamma$ on $V$ resp. $V_0$.*

*Then $E_0$ is a holomorphically trivial vector bundle and $\Gamma(X, E_0) = \Gamma(X, E)$.*

For antiholomorphic representations we are now able to give a precise description of the global sections of the associated vector bundle.

**Proposition 9.** *Let $G$ be a connected complex Lie group, $\Gamma$ a lattice, $V$ a complex vector space and $\rho : G \to GL(V)$ be an antiholomorphic representation. Let $E$ denote the flat vector bundle over $G/\Gamma$ which is induced by $\rho|_\Gamma$.*

*Let $\Sigma$ denote the set of all vectors $v \in V$ which are invariant under $\rho(G')$ and a $\rho(\gamma)$-eigenvector with real eigenvalue for every $\gamma \in \Gamma$.*

*Then $H^0(G/\Gamma, E) \simeq <\Sigma>_\mathbb{C}$ where $<\Sigma>_\mathbb{C}$ denotes the complex vector space spanned by $\Sigma$.*

*Proof.* By construction $<\Sigma>_\mathbb{C}$ is the largest $\rho(G)$-invariant vector subspace of $V$ inducing a holomorphically trivial vector subbundle of $E$. □

## 11. Subbundles of flat bundles

**Theorem 6.** *Let $G$ be a connected complex Lie group, $\Gamma$ a lattice and $E \to G/\Gamma$ be a flat vector bundle given by an essentially antiholomorphic representation $\rho : G \to GL(n, \mathbb{C})$.*

*Assume that every meromorphic function on $X = G/\Gamma$ is constant.*

*Let $L \subset E$ be a vector subbundle. Then $L$ also admits a flat connection.*

*If $G = G'$, then $L$ is parallel with respect to the given flat connection on $E$.*

*Proof.* Let $k = rank(L)$. By passing to $\Lambda^k L \subset \Lambda^k E$ we may assume that $L$ is a line bundle. The flat connection on $E$ yields a canonical way to lift the $G$-action on $X$ to a $G$-action on $E$. The subbundle $E_0 = \oplus_{g \in G} g^* L$ is a $G$-invariant subbundle of $E$ and therefore parallel with respect to the flat connection on $E$. Moreover it is given by an essentially antiholomorphic representation (see cor. 3). We may thus assume $E_0 = E$. Let $d = rank(E_0)$ and choose $g_0, \ldots, g_d \in G$ such that the line subbundles $L_i = g_i^* L$ are in general position, i.e. such that for every number $m$ with $0 \leq m \leq d$ and every choice of $0 \leq i_1 < \ldots i_m \leq d$ the subsheaf of $E$ spanned by $L_{i_1} + \ldots + L_{i_m}$ has rank $m$.

**Claim 9.** *The line bundles $L_0, \ldots, l_d$ are in general position at every point and therefore yield a trivialization of $\mathbb{P}(E_0)$.*

*Proof.* Let $U$ be an open neighbourhood of a point such that all the $L_i$ admit nowhere vanishing sections $\sigma_i$. For every $k \in \{0, \ldots, d\}$ we obtain a section $\alpha_k$ in $\wedge^d E_0$ by $\alpha_k = \wedge_{i \neq k} \sigma_i$. Let $\alpha = \otimes_k \alpha_k$. Then $\alpha$ is a section in the line bundle $\otimes^{d+1}(\wedge^d E_0)$ and is vanishing exactly where the $L_0, \ldots, L_d$ fail to be in general position. Thus the set of all points in $X$ where the $L_0, \ldots, L_d$ are not in general position constitute an analytic hypersurface of $X$. But the assumptions on $X$ imply that $X$ contains no hypersurfaces. Hence $L_0, \ldots, L_d$ are in general position everywhere. $\square$

Thus $\mathbb{P}(E_0)$ is a holomorphically trivial $\mathbb{P}_{d-1}$-bundle. Observe that $\mathbb{P}(E_0)$ is defined by an essentially antiholomorphic group homomorphism $\tilde{\rho}: \Gamma \to PGL(d, \mathbb{C})$. Using thm. 2 it follows that $\tilde{\rho}|_{G' \cap \Gamma}$ is trivial. Hence $\tilde{\rho}(\Gamma)$ constitutes a commutative subgroup of $PGL(d, \mathbb{C})$. This implies that there is a fixed point for the $\tilde{\rho}(\Gamma)$-action on $\mathbb{P}(V)$. Corresponding to this fixed point there exists a sub-line bundle $L' \subset E$ which is $G$-invariant and therefore parallel to the given flat connection on $E$.

**Claim 10.** *The line bundles $L'$ and $L$ are isomorphic (as holomorphic line bundles).*

*Proof.* Recall that $\mathbb{P}(E)$ is holomorphically trivial and that every meromorphic function on $X$ is constant. It follows that there is a unique trivialization of $\mathbb{P}(E)$ and that every section of $\mathbb{P}(E) \to X$ is constant with respect to this trivialization. This implies that any two line subbundles of $E$ are holomorphically isomorphic. $\square$

Thus $L$ is isomorphic to a flat line bundle. However, it is parallel with respect to the given flat connection on $E$ if and only if it is $G$-invariant. This in turn is equivalent to the assertion that $L$ corresponds to a fixed point of the $\Gamma$-action on $\mathbb{P}(V)$. This action is trivial, if $G = G'$. $\square$

JÖRG WINKELMANN, MATH. INSTITUT NA 4/69, RUHR-UNIVERSITÄT BOCHUM, 44780 BOCHUM, GERMANY
  *E-mail address*: jw@cplx.ruhr-uni-bochum.de